\documentclass[journal]{IEEEtran}

\usepackage{cite}
\usepackage{amsmath}
\usepackage{url}
\usepackage{graphicx}
\usepackage{color}
\usepackage{steinmetz}
\usepackage{placeins}
\usepackage{float}
\usepackage{citesort}
\usepackage{tabularx,colortbl}
\usepackage{ifthen}

\makeatletter

\newcounter{author}
\renewcommand{\author}[2][]{
   \stepcounter{author}
   \@namedef{author@\theauthor}{#2}
   \@namedef{authorlabel@\theauthor}{#1}
}

\newcounter{address}
\newcommand{\address}[2][]{
   \stepcounter{address}
   \@namedef{address@\theaddress}{#2}
   \@namedef{addresslabel@\theaddress}{#1}
}

\newcommand{\alsep}{and}

\def\newmaketitle{\par%
  \begingroup%
  \normalfont%
  \def\thefootnote{}
  \def\footnotemark{}
  \let\@makefnmark\relax
  \footnotesize
  \footnotesep 0.7\baselineskip
  \normalsize%
  \twocolumn[\thenewmaketitle\@IEEEaftertitletext]%
  \if@IEEEusingpubid
     \enlargethispage{-\@IEEEpubidpullup}%
  \fi
  \endgroup
  \setcounter{footnote}{0}\let\maketitle\relax\let\@maketitle\relax
  \gdef\@thanks{}%
  \let\thanks\relax}

\def\thenewmaketitle{
  \newpage
  \begin{center}%
    \vskip0.2em{\Huge\@IEEEcompsoconly{\sffamily}\@IEEEcompsocconfonly{\normalfont\normalsize\vskip 2\@IEEEnormalsizeunitybaselineskip
   \bfseries\large}\@title\par}\vskip1.0em\par%
    \vspace{1ex}
    \newcounter{c@author}
    \newcounter{c@tmp}
    \ifthenelse{\value{author}=2}{%
      \newcommand{\liand}{ and }}{%
      \newcommand{\liand}{, and }}
    \ifthenelse{\value{address}<2}{%
      \@nameuse{author@1}%
      \stepcounter{c@author}%
      \whiledo{\value{c@author}<\value{author}}{%
        \setcounter{c@tmp}{\value{author}}%
        \addtocounter{c@tmp}{-\value{c@author}}%
        \ifthenelse{\value{c@tmp}=1}{%
          \renewcommand{\alsep}{\liand}}{\renewcommand{\alsep}{, }}%
        \stepcounter{c@author}\alsep \@nameuse{author@\thec@author}}\\%
    }
    {
      \@nameuse{author@1}${}^{(\ref{\@nameuse{authorlabel@1}})}$%
      \stepcounter{c@author}%
      \whiledo{\value{c@author}<\value{author}}{%
      \setcounter{c@tmp}{\value{author}}%
      \addtocounter{c@tmp}{-\value{c@author}}%
      \ifthenelse{\value{c@tmp}=1}{%
        \renewcommand{\alsep}{\liand}}{\renewcommand{\alsep}{, }}%
      \stepcounter{c@author}\alsep \@nameuse{author@\thec@author}%
        ${}^{(\ref{\@nameuse{authorlabel@\thec@author}})}$%
      }
    }
    \vspace{0.2ex}

    \ifthenelse{\value{address}>0}{%
      \ifthenelse{\value{address}=1}{
        {\@nameuse{address@1}}
      }
      {
        \newcounter{c@address}

        \begin{center}
        \whiledo{\value{c@address}<\value{address}}
        {
          \refstepcounter{c@address}
            ${}^{(\thec@address)}$\,%
              \label{\@nameuse{addresslabel@\thec@address}}%
              \@nameuse{address@\thec@address}\\ %
        }
        \end{center}
      } 
    }
    {
      \relax
    }
  \end{center}
}

\makeatother

\title{A Transient Port-Extraction Technique for Antenna Feed Optimization}

\author[org1]{Sean F. DePalma}
\author[org1]{Omkar H. Ramachandran}
\author[org1]{Leo C. Kempel}
\author[org2]{B. Shanker}

\address[org1]{Department of Electrical and Computer Engineering, Michigan State University, East Lansing, MI 48824}
\address[org2]{Department of Electrical and Computer Engineering, The Ohio State University, Columbus, OH 43210}

\begin{document}

\newmaketitle

\begin{abstract}

Optimization of strongly non-linear tightly coupled feeds attached to antennas is a challenging problem from a purely computational perspective. One can imagine that an optimization would (a) need to be in the time domain, and (b) has to be self-consistently coupled with the linear antenna (or electromagnetic) system. These two imply that the cost of optimization is governed by the need to repeatedly evaluate the fully coupled cost function. This paper leverages a recently developed transient port-extraction technique to circumvent this challenge and is \emph{agnostic} to the optimization scheme. This approach provides a representation of the entire linear electromagnetic system at the port and can readily integrate with \emph{any} non-linear circuit analysis and optimization methodology. In this paper, we
demonstrate optimization of linear and non-linear circuit feed parameters that are tightly coupled to broadband radiating systems.

\end{abstract}

\section{Introduction}
The capability to efficiently optimize a radio frequency (RF) network consisting of a broadband electromagnetic (EM) device fed by strongly non-linear lumped circuit components is immediately relevant in modern RF design.
Such a method would allow for assessment and mitigation of radiative coupling effects early in the design of high frequency devices with compact features.
Due to the strongly nonlinear nature of the devices of interest, we focus primarily on transient methods.

The current state of the art to analyse broadband devices excited by lumped circuit sources \cite{jin_riley,jin2015finite} relies primarily on a tight coupling between the linear, broadband EM system interacting with a set of lumped nonlinear devices.
While such methods are accurate against measured data, the computational cost of evaluating nonlinear iterations can become prohibitively expensive.
More recently, transient parameter extraction methods have been proposed \cite{scott_pe,Ramachandran} that greatly amortize this cost while retaining the accuracy of a tightly coupled solution.
Furthermore, the port responses were shown to be agnostic to circuit model, allowing for reuse with different configurations. In characterizing the EM system in this way, the reduced-order EM-circuit network is amenable to efficient non-linear feed optimization following \emph{any} method so long as the EM system remains unchanged. Problems related to non-linear optimization of the antenna feeding circuit are actively studied in the literature, primarily for filtering, beam-steering, reflection coefficient improvement, and other radiative characteristics of complex arrays \cite{haupt}.

In this letter, we use transient parameter extraction to construct a multi-objective optimization method for nonlinear EM devices. We validate our method through a number of numerical examples, before using a genetic algorithm to reduce distortion from a strongly non-linear amplifying circuit feeding a broadband antenna. Likewise, we demonstrate signal to interference and noise ratio (SINR) improvement by implementing a reverse beam-steering scheme to mitigate time-varying interference exclusively through circuit feed modification.
 
\section{Problem Statement}
A radiating object contains $N_{p}$ ports each associated with a lumped circuit subsystem. A transient response $G_{i}(t)$ is extracted with each port $i$, following the method detailed in \cite{Ramachandran}. 
Each attached circuit subsystem, represented is assumed to be ideal and tunable.
Our goal is to modify the circuit subsystems feeding the unchanged EM system to optimize a cost function derived from the overall EM-circuit response. 

The tunable components are represented by a vector of design variables $\mathbf{\mathcal{E}}$. 
These are optimized with the goal of minimizing a cost function $J(\mathbf{\mathcal{E}})$ defined as 
\begin{equation}
    J(\mathbf{\mathcal{E}})= |{d - f(\mathbf{\mathcal{E}})}|,
    \label{optfunc}
\end{equation}
where $f(\mathbf{\mathcal{E}})$ is a nonlinear function dependent on the general EM-circuit response and $d$ represents a desired system response.

\section{Methodology} 

Let $\mathcal{L}$ and $\mathcal{C_{EM}}$ refer respectively to an abstraction of a linear Maxwell operator and the coupling relation between this operator and a lumped port feed model. For the results in this work, $\mathcal{L}$ refers to a finite element discretization of the EM system. The lumped circuit subsystems are analysed using a spice-like nodal network \cite{mna}, and both systems is solved temporally using an Newmark-$\beta$ scheme \cite{newmark}, with $\gamma=0.5$ and $\beta=0.25$, a configuration known to be unconditionally stable. The time varying field coefficients in the EM system $\bar{E}$ relate to the attached circuit currents $\bar{J}$ through $\mathcal{L}\circ\left[\bar{E}(t)\right] + \mathcal{C_{EM}}\circ\left[\bar{J}(t)\right] = 0$.
We represent $\bar{E}(t)$ and $\bar{J}(t)$ using the same set of temporal basis functions $\lambda(t)$: $\bar{E}(t) = \sum_{i=0}^{N_{ts}} e_{i}\lambda(t-t_{i}) \; \text{and} \; \bar{J}(t) = \sum_{i=0}^{N_{ts}} j_{i}\lambda(t-t_{i})$.
The port response $G(t)$ can be found due to a single temporal basis function $\lambda(t-t_{i})$: $\mathcal{L}\circ G(t) = -\mathcal{C_{EM}} \circ \lambda (t)$.
This response can then be discretely convolved with \emph{any} circuit current $\bar{J}$ to obtain the resulting fields, i.e $\bar{E}(t)=\bar{J}(t)\star\bar{G}(t)$. Complete details of this method can be found in \cite{Ramachandran}. We note that in systems where the EM system has significantly more degrees of freedom than the lumped circuit system, constraining nonlinear iterations to the circuit results in significant computational speedup.

Consider an antenna array with incident desired, interfering, and noise signals. Our objective is to maximize the signal to interference $\&$ noise ratio (SINR) using Applebaum's criterion \cite{Applebaum}. This has been applied within various analytic optimization cases, including the use of a genetic algorithm (GA) to vary digital phase shift in an array of isotropic sources \cite{Weile}. Our objective is to generalize this to a finite-element representation of a physical antenna.

For a linear array with $N$ equally spaced identical elements, the source amplitudes and phases are defined with a weighting vector $\mathbf{w} = \left|{w_{m}}e^{j\phi_{m}}\right|$
where $0\leq m \leq N$, and $w_{m}$ and $\phi_{m}$ are the magnitude and phase respectively for each element $m$. We define an incident angle vector: $\mathbf{u}(\theta)  = \left[ \; 1, \; e^{j\pi sin(\theta)}, \; \dots, \; e^{(N-1)j\pi sin(\theta)} \; \right] $.

A desired signal is incident upon the array along with time-changing interference and element-wise noise, combining as a total received signal 
$\mathbf{x}(t) = \mathbf{x}_{\text{desired}}(t) + \mathbf{x}_{\text{interference}}(t) + \mathbf{x}_{\text{noise}}(t)$.
The frequency distribution of desired signal and its angle of incidence $\theta_{d}$ are assumed to be known. This signal is represented as $    \mathbf{x}_{\text{desired}}(t) = a_{d}(t) \mathbf{u}(\theta_{d})$.
where $a_{d}(t)$ is a known time-varying desired signal. The interference $\mathbf{x}_{\text{i}}(t)$ is defined in the same form as $\mathbf{x}_{\text{desired}}(t)$ with $a_{i}(t)$ time-varying interference signals with varied incident angles $\theta_{i}$. The element-wise noise $\mathbf{x}_{\text{n}}(t)$ is introduced to the system as incident Gaussian noise. Continuing to define quantities as done in \cite{Weile}, a generalized covariance matrix is defined as $\Phi_{q} = E\{ \mathbf{x}_{q}(t)^{*} \, \mathbf{x}_{q}(t)^{T} \}$.
in which $q$ represents one of the desired, interference, or noise components and $E \{ \cdot \}$ is the expectation operator. Signal power is calculated as $P_{q} = \mathbf{w}^{\mathcal{H}} \Phi_{q} \, \mathbf{w}$.
with the superscript $\mathcal{H}$ denoting the Hermitian operator.

The SINR is calculated as the ratio between the power received from the desired signal against that of the total interference plus noise. Consider the desired signal has a known amplitude $a$ used to define $\bar{a}^{2} = E \{ \left| a \right| ^{2} \}$. The Applebaum criterion from \cite{Applebaum} then maximizes SINR from the defined quantities:
\begin{equation}
    \label{array_opti_costfunc}
    \text{SINR} = \frac{P_{\text{desired}}}{P_{\text{undesired}}}
    = \bar{a}^{2} \frac{\left| \mathbf{w}^{T} \mathbf{u}(\theta_{d})^{2} \right|}{\mathbf{w}^{\mathcal{H}}(\Phi_{i}+\Phi_{n}) \, \mathbf{w}}
\end{equation}
The incident angles $\theta_{d}$ of the desired is known but the amplitudes of the signal and interference become inseparable upon incidence, meaning it is not possible to directly compute $\bar{a}^{2}$ and $\Phi_{i} + \Phi_{n}$ from the total received signal $\mathbf{x}(t)$. 
Thus, we define $J(\mathbf{\mathcal{E}})$ for this problem based on a simplified function which has the same maximum as \eqref{array_opti_costfunc} when using the weighting vector $\mathbf{w}$ as the design vector $\mathbf{\mathcal{E}}$ by omitting $\bar{a}^{2}$ and adding $\Phi_{d}$ to $\Phi_{i} + \Phi_{n}$ in the denominator \cite{Weile}.
With these definitions in place for the antenna array, we will optimize the array system for SINR maximization.

\section{Results \label{sec:results}}
For the results presented in this section, the voltage sources used are defined by $v(t)=\cos(2\pi f_{0}t)\exp\left(-(t-6\sigma)^2/2\sigma^2 \right)$, where $\sigma=3/(2\pi f_{\text{bw}})$ for center frequency $f_{0}$ and bandwidth $f_{\text{bw}}$. Likewise, the stepsize is chosen as $\Delta_{t}=(20(f_{0}+f_{\text{bw}}))^{-1}$. 
In each case, the simulation domain is truncated by a PML placed $\frac{1}{2} \lambda_{max}$ from each side of the device. The resulting geometry is discretized using tetrahedra with an average side length of $\lambda_{min}/30$.
After the EM system is port-extracted, circuit components will be modified and cost function minimized using a GA from the pymoo python module \cite{pymoo}.

\subsection{Comparison Using Linear System} 
\label{section:LinearAnalysisVivaldi}
\begin{figure}
    \centering
    \includegraphics[width=0.48\textwidth]{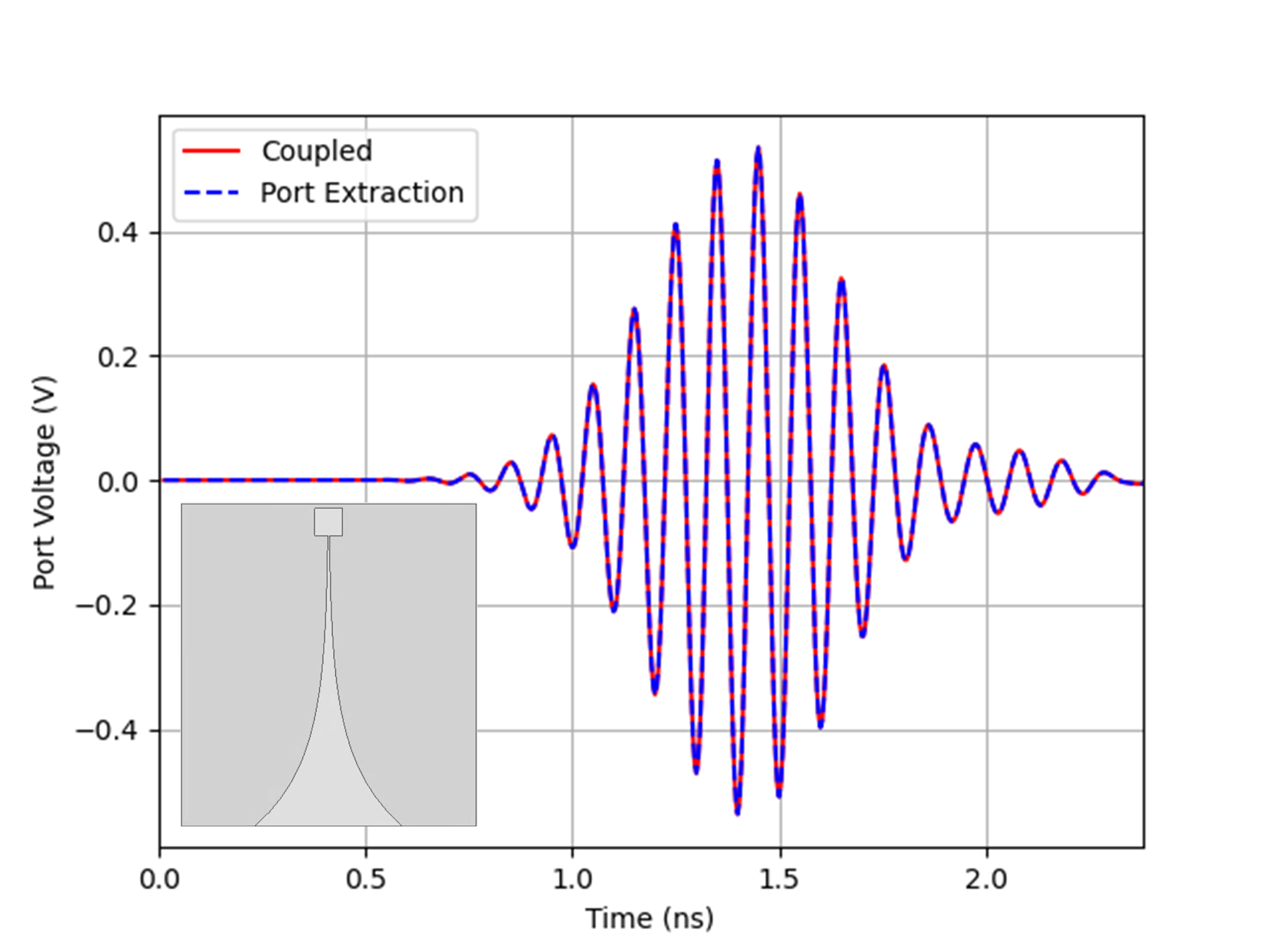}
    \caption{Voltage at port nodes of an antenna found with self-consistent analysis and transient port-extraction.}
    \label{fig:vivaldi_VoltagePEvsCoupled}
\end{figure}

A Vivaldi antenna is shown as an inset in Fig. \ref{fig:vivaldi_VoltagePEvsCoupled}. It is fed by a linear circuit and this system is simulated to compare the results of transient port-extraction and self-consistent analysis.
The face of the antenna consists of symmetrical PEC exponential curves over a lossless 
substrate with relative permittivity $\epsilon_{r}$ = 2.2. The substrate has side lengths 41mm and height 1.575mm. 
The feeding circuit is connected to the antenna across a single edge located at the base of the exponential flare, running from halfway through the substrate to one PEC face of the antenna. This circuit contains a voltage source with $f_{0}$ = 10GHz and  $f_{\text{bw}}$ = 2GHz with a series 50$\Omega$ resistor. 
The system is analyzed using port-extraction and self-consistent methods with an iterative solver tolerance of $10^{-12}$. The voltage at the port is measured for each method and compared in Fig. \ref{fig:vivaldi_VoltagePEvsCoupled}. The $\mathcal{L}^{2}$ error between the port voltage of the two methods is $6.4\times10^{-12}$, demonstrating solver-precision agreement between transient port-extraction and self-consistent analysis.

\subsection{Multiple-Objective Optimization of Linear System}
\begin{figure}
    \centering
    \includegraphics[width=.9
    \linewidth]{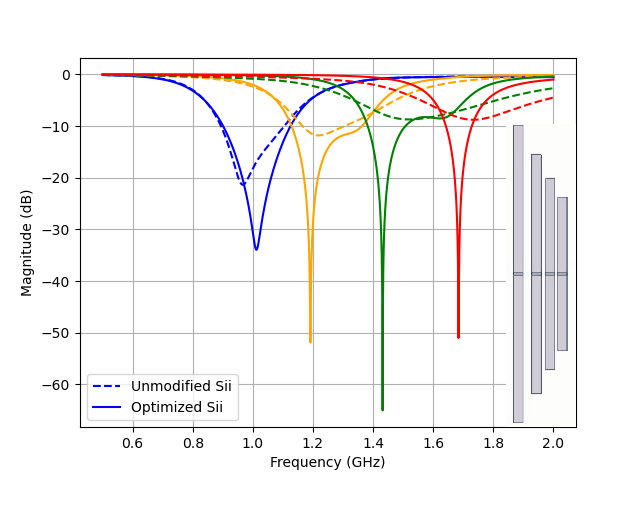}
\caption{Optimized Reflection Coefficients of Log-Periodic Dipole Array}
\label{plot:logpdipolesOptimized}
\end{figure}

A four-element log periodic dipole array is simulated to demonstrate multiple-objective optimization through individually minimizing reflection coefficients. The geometry is shown as an inset in Fig.  \ref{plot:logpdipolesOptimized}. The first dipole is defined with a total length of 14.9cm, the second with 11.9cm, the third with 9.5cm, and the fourth with 7.6cm. The center of the second dipole is 0.9cm away from the center of the first dipole, the third is 1.6cm away from the first, and the fourth is 2.2cm away from the first. Each dipole has a diameter of 4mm. 
A vertical gap of 1mm is added between the flanges of each dipole. A circuit subsystem is connected across each gap for each dipole. Each circuit consists of a voltage source with a $f_{0}$ = 1.25GHz and $f_{\text{bw}}$ = 750MHz. A resistor with an initial value of 50$\Omega$ is connected in series with each voltage source. The geometry is discretized using a mesh with 1.8M elements. 

The reflection coefficients of each dipole are minimized using the GA to modify resistor values. The settings modified in the default algorithm are an initial population of 80, 45 offspring per generation, and 75 generations. Using \eqref{optfunc} as the cost function, $J(\mathbf{\mathcal{E}})$ is minimized with $f(\mathbf{\mathcal{E}})$ defined as the vector of measured reflection coefficients with $d=0$. The input vector $\mathbf{\mathcal{E}}$ contains the resistances in each source circuit. From a set of solutions produced, the optimal solution was chosen as the minimum sum of the individual reflection coefficients. This result is compared to $50\ \Omega$ feeds in Fig. \ref{plot:logpdipolesOptimized}.

\subsection{Non-linear Optimization of a Broadband System}

\begin{figure}
    \centering
    \includegraphics[width=.9
    \linewidth]{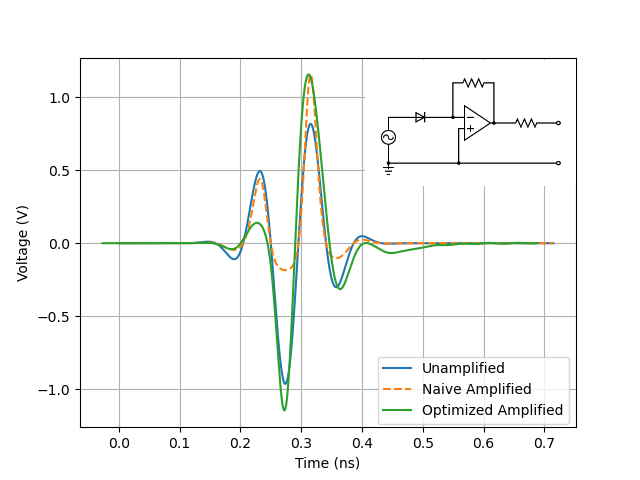}
\caption{Initial and optimized amplified output compared to unamplified in time-domain. Inset is naive non-linear amplifier circuit.}
\label{plot:NLBroadbandResult2}
\end{figure}

\begin{figure}
    \centering
    \includegraphics[width=.9
    \linewidth]{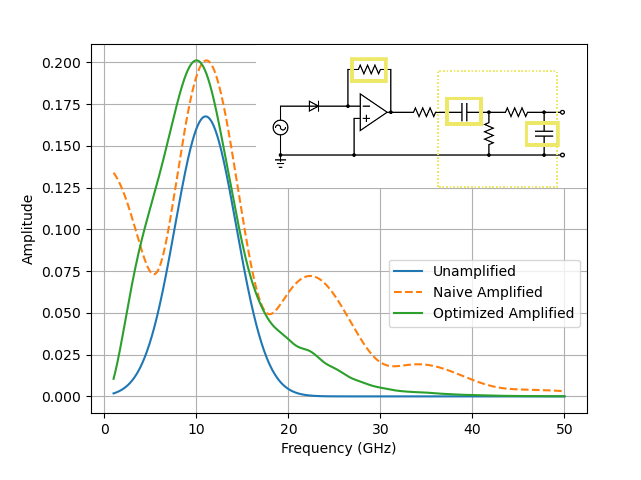}
\caption{Initial and optimized amplified output compared to unamplified in frequency-domain. Inset is optimized non-linear amplifier circuit with modified components highlighted.}
\label{plot:NLBroadbandResult3}
\end{figure}

\begin{figure}
    \centering
    \includegraphics[width=.7
    \linewidth]{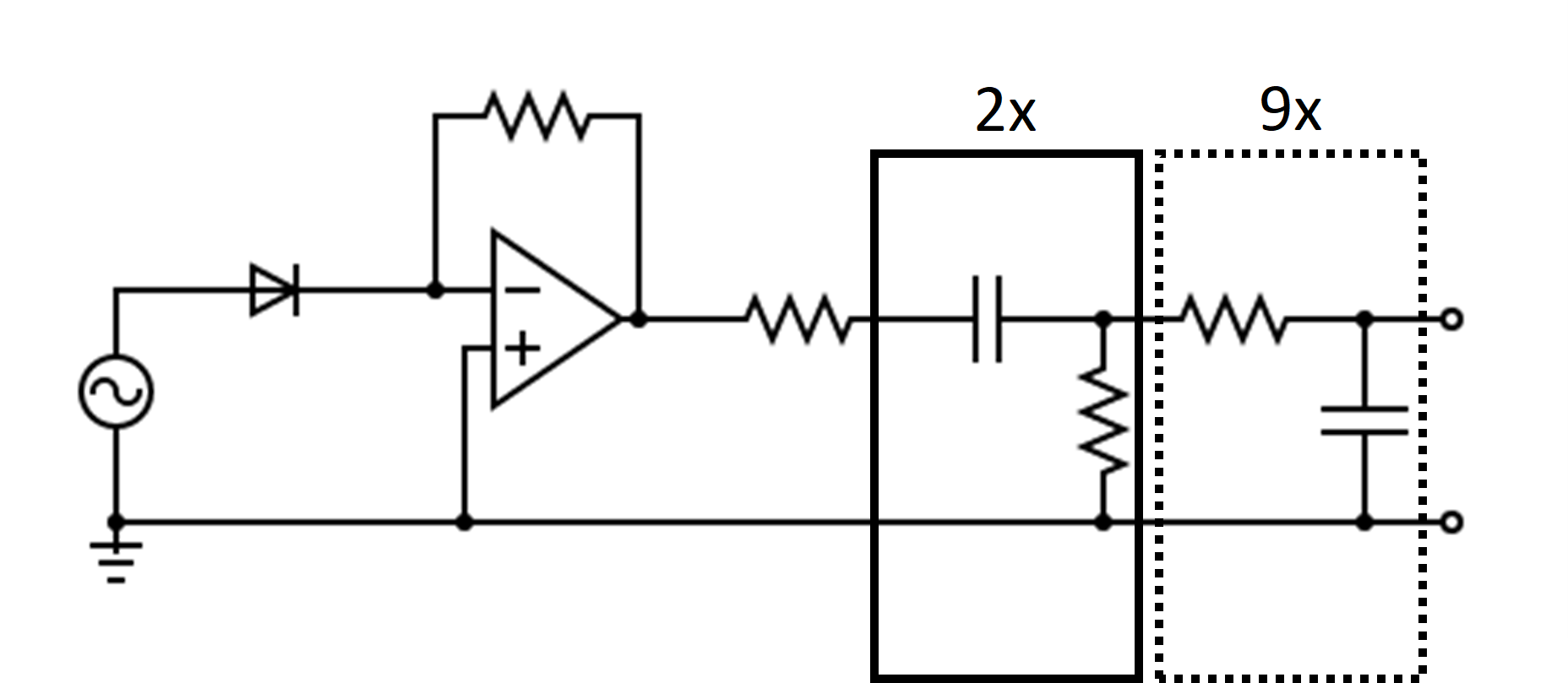}
\caption{Resulting optimized circuit with the solid box representing nine repeated RC low-pass filters and the dashed box representing two repeated RC high-pass filters}
\label{plot:NLBroadbandResult4}
\end{figure}

The Vivaldi antenna investigated within section 
\ref{section:LinearAnalysisVivaldi}
is now used for non-linear broadband analysis. The port response from that example was reused with the circuit topology modified to the exponential amplifier inset in Fig. \ref{plot:NLBroadbandResult2}. The voltage source is defined for $f_{0}$ = 11GHz and $f_{\text{bw}}$ = 10GHz.
The ideal diode in the circuit is modeled using characteristics from the 1N6263 Schottky diode. The initial amplifying resistance is set as 1M$\Omega$ and the output resistance as 100$\Omega$. Feeding the Vivaldi antenna with this circuit produces the "naive", or non-optimized, amplified time response shown in Fig. \ref{plot:NLBroadbandResult2} and likewise frequency response in Fig. \ref{plot:NLBroadbandResult3}, both of which are compared to the result of feeding the antenna with the unamplified voltage source and series 50$\Omega$ resistor. As is evident, the dipole excites higher harmonics, which we seek to suppress by means of filtering the output. A passive band-pass filter is introduced to do this. The resulting circuit configuration is shown inset in Fig. \ref{plot:NLBroadbandResult3}. A cost function $J(\mathbf{\mathcal{E}})$ is determined by minimizing the relative $\mathcal{L}^{1}$ error of the norm of the Fourier-transformed port voltage, minimizing the difference in frequency-domain between the unamplified response and the filtered amplified response. The design space $\mathbf{\mathcal{E}}$ consists of varying the amplifying resistance from 100$\Omega$ to 1M$\Omega$, both capacitors individually from 1pF to 100pF, and the number of low-pass and high-pass RC filters from one to ten. Setting the initial population to 150, offspring to 50, and generations to 30 for this optimization scheme produced the results denoted as optimized amplified shown in Fig. \ref{plot:NLBroadbandResult2} and \ref{plot:NLBroadbandResult3}. This resulted in an optimal circuit topology consisting of two identical high-pass RC pairs of 1$\Omega$ and 5.4pF and nine identical two-pass RC pairs of 1$\Omega$ and 6.7pF and is represented in Fig. \ref{plot:NLBroadbandResult4}. 
The normalized relative $\mathcal{L}^{1}$ error in the frequency spread is reduced from 0.76 for naive amplification compared to 0.16 produced by the proposed method.

\subsection{Adaptive SINR Optimization by Reverse Beam-Forming}
The circuit system of a five-element linear dipole system is optimized next as an adaptive array. This stems from a previous investigation into maximizing SINR in a system with incident desired and interfering waves using the Applebaum criterion \cite{Weile}. It utilized a linear array of isotropic sources with digital phase shifts and optimization not involving the simulation of an EM system. The work presented here uses a linear dipole array with phase shift induced by analog circuit components. The phase shift is measured using a time difference between the voltage source and port voltages found after port-extracting the EM system and feeding each dipole with a non-linear circuit. 

The EM system consists of five dipoles with center frequency $f_{0}$ = 1.375GHz with spacing and length equal to one half-wavelength, and a bandwidth $f_{\text{bw}}$ = 100MHz. The diameter of the dipoles is 4mm and a 1mm gap is added between each set of flanges over which each source spans. 

\begin{figure}
    \centering
    \includegraphics[width=0.45\textwidth]{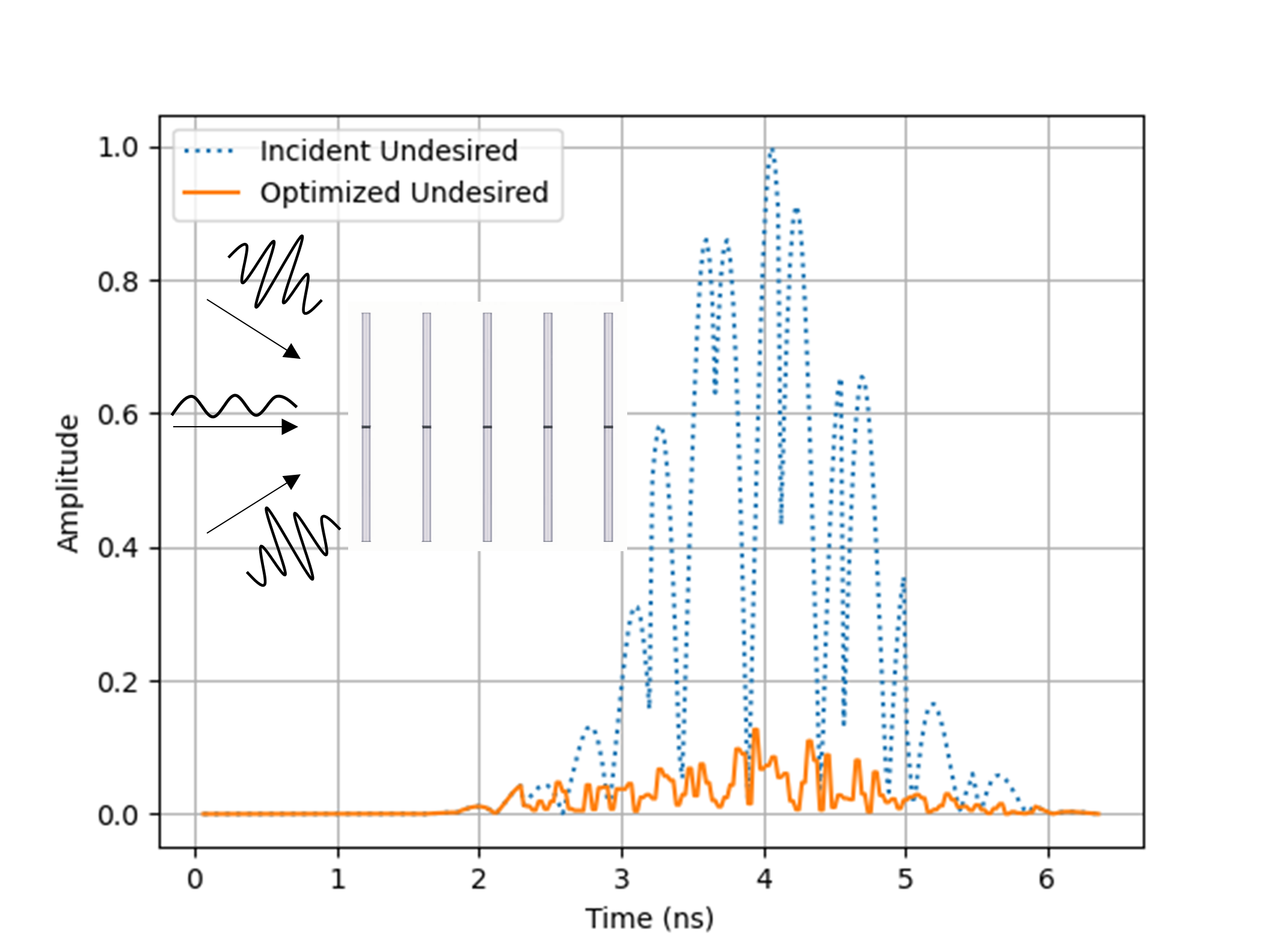}
    \caption{Comparing the received and mitigated amplitudes of an undesired signal by a five-element dipole array. Inset depicts the array with incident signals.}
    \label{fig:SINR_Cancellation}
\end{figure}

The amplitude of the interference is defined as 20dB greater than that of the desired signal, and the element-wise noise is defined as 20dB less than the desired signal. We define a single interfering wave incident upon this five-element array, with the angle of the incidence changing between $40^{\circ}$ and $-30^{\circ}$ off broadside every 50 time steps. The desired signal is incident broadside of the array. Both signals are defined as modulated Gaussians. 
The voltage source amplitudes in the linear dipole array are fixed to produce a Dolph-Chebyshev beam pattern broadside of the array with 20dB sidelobe suppression. 
The SINR is optimized using a cost function $J(\mathbf{\mathcal{E}})$ equal to the modified version of \eqref{array_opti_costfunc}. 

\begin{figure}
    \centering
    \includegraphics[width=0.45\textwidth]{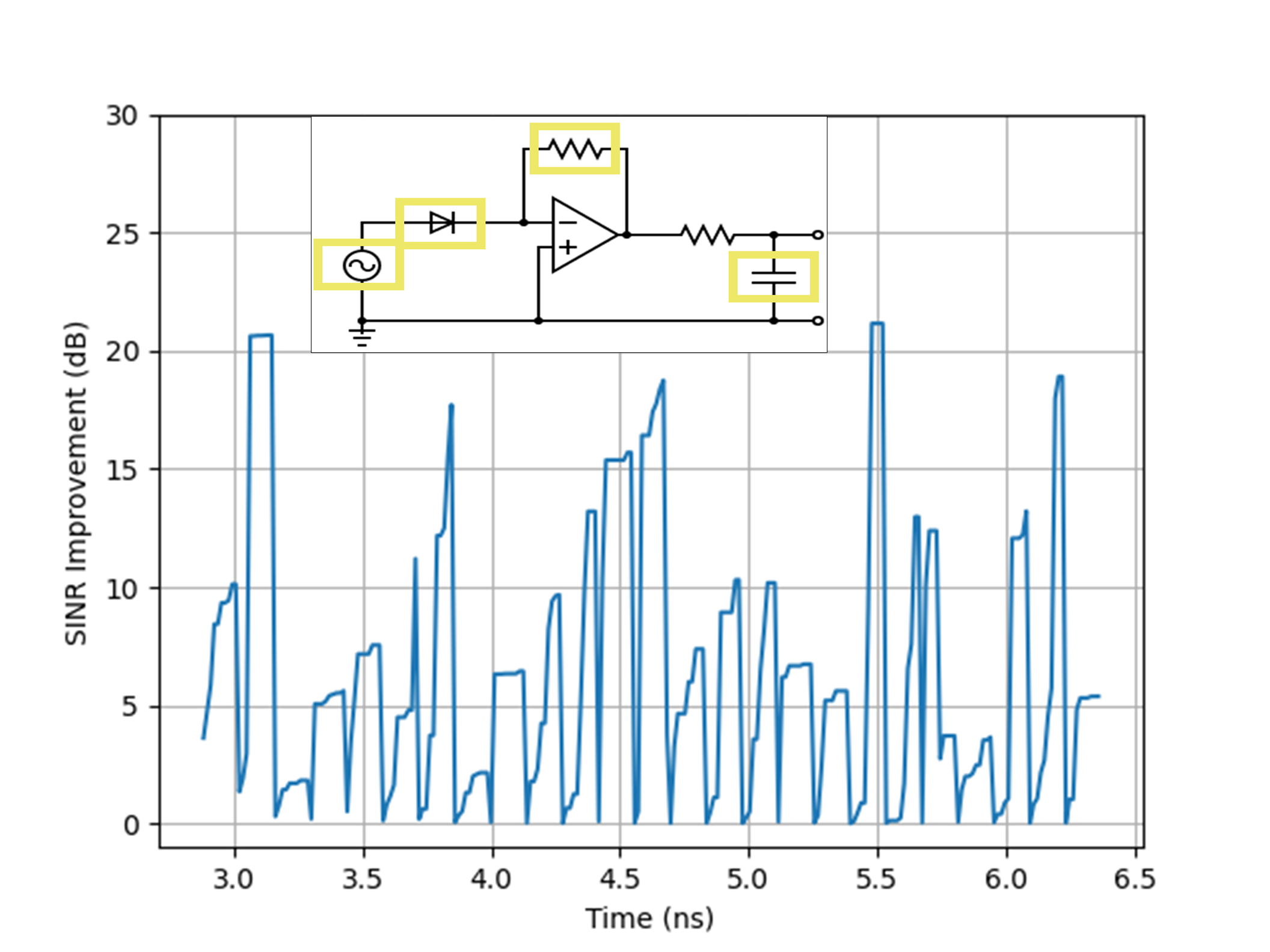}
    \caption{Improvement in the SINR as a function of time. The inset Fig. depicts the non-linear feeding circuit, with tunable components highlighted.}
    \label{fig:SINR_NL}
\end{figure}

This scheme is applied with strongly non-linear circuits feeding each dipole. The circuit stencil used is inset in Fig. \ref{fig:SINR_NL}, with tunable components highlighted. The design space consists of a voltage source delay between 0 to 150 time steps, an amplifying resistor varied from 100$\Omega$ to 1M$\Omega$, a filtering capacitor varied from 30pF to 80pF, and diode thermal voltage modified by varying temperature between 285 Kelvin and 300 Kelvin. The resulting SINR improvement over time is shown in Fig. \ref{fig:SINR_NL}. This optimization scheme using temporal port extraction demonstrates a computational cost improvement of approximately five orders of magnitude, improving from 4.5s for a single self-consistent forward solve of the system to a post-extraction circuit solve completing in 243$\mu$s.

\section{Conclusion}
In this letter, we demonstrate the use of a recently developed transient port extraction method to optimize broadband antennas tightly coupled to strongly non-linear feeding circuits. We demonstrate that the technique is amenable to non linear optimization through a variety of numerical examples utilizing a genetic algorithm to modify parameters in feeding circuits. We demonstrate significant speedup when using port-extraction, as the optimization is constrained to the circuit solve. 

\bibliographystyle{IEEEtran}
\bibliography{references}

\begin{thebibliography}{10}
\providecommand{\url}[1]{#1}
\csname url@samestyle\endcsname
\providecommand{\newblock}{\relax}
\providecommand{\bibinfo}[2]{#2}
\providecommand{\BIBentrySTDinterwordspacing}{\spaceskip=0pt\relax}
\providecommand{\BIBentryALTinterwordstretchfactor}{4}
\providecommand{\BIBentryALTinterwordspacing}{\spaceskip=\fontdimen2\font plus
\BIBentryALTinterwordstretchfactor\fontdimen3\font minus
  \fontdimen4\font\relax}
\providecommand{\BIBforeignlanguage}[2]{{%
\expandafter\ifx\csname l@#1\endcsname\relax
\typeout{** WARNING: IEEEtran.bst: No hyphenation pattern has been}%
\typeout{** loaded for the language `#1'. Using the pattern for}%
\typeout{** the default language instead.}%
\else
\language=\csname l@#1\endcsname
\fi
#2}}
\providecommand{\BIBdecl}{\relax}
\BIBdecl

\bibitem{jin_riley}
J.~Jin and D.~J. Riley, \emph{Finite Element Analysis of Antennas and
  Arrays}.\hskip 1em plus 0.5em minus 0.4em\relax Wiley-IEEE Press, 2009.

\bibitem{jin2015finite}
J.-M. Jin, \emph{The finite element method in electromagnetics}.\hskip 1em plus
  0.5em minus 0.4em\relax John Wiley \& Sons, 2015.

\bibitem{scott_pe}
S.~{O’Connor}, S.~{Hughey}, D.~{Dault}, A.~J. {Pray}, J.~M. {Villa-Giron},
  and B.~{Shanker}, ``A novel port/network parameter extraction technique for
  coupling circuits with full-wave time-domain integral equation solvers,''
  \emph{IEEE Transactions on Microwave Theory and Techniques}, vol.~67, no.~2,
  pp. 553--564, 2019.

\bibitem{Ramachandran}
O.~H. Ramachandran, S.~O’Connor, Z.~D. Crawford, L.~C. Kempel, and
  B.~Shanker, ``Port parameter extraction-based self-consistent coupled
  em-circuit fem solvers,'' \emph{IEEE Transactions on Components, Packaging
  and Manufacturing Technology}, vol.~12, no.~6, pp. 1040--1048, 2022.

\bibitem{haupt}
R.~Haupt and D.~Werner, \emph{Genetic Algorithms in Electromagnetics}.\hskip
  1em plus 0.5em minus 0.4em\relax John Wiley \& Sons, Ltd, 2007.

\bibitem{mna}
{Chung-Wen Ho}, A.~{Ruehli}, and P.~{Brennan}, ``The modified nodal approach to
  network analysis,'' \emph{IEEE Transactions on Circuits and Systems},
  vol.~22, no.~6, pp. 504--509, 1975.

\bibitem{newmark}
O.~C. Zienkiewicz, ``A new look at the newmark, houbolt and other time stepping
  formulas. a weighted residual approach,'' \emph{Earthquake Engineering \&
  Structural Dynamics}, vol.~5, no.~4, pp. 413--418, 1977.

\bibitem{Applebaum}
S.~Applebaum, ``Adaptive arrays,'' \emph{IEEE Transactions on Antennas and
  Propagation}, vol.~24, no.~5, pp. 585--598, 1976.

\bibitem{Weile}
D.~Weile and E.~Michielssen, ``The control of adaptive antenna arrays with
  genetic algorithms using dominance and diploidy,'' \emph{IEEE Transactions on
  Antennas and Propagation}, vol.~49, no.~10, pp. 1424--1433, 2001.

\bibitem{pymoo}
J.~{Blank} and K.~{Deb}, ``pymoo: Multi-objective optimization in python,''
  \emph{IEEE Access}, vol.~8, pp. 89\,497--89\,509, 2020.

\end{thebibliography}



%

\end{document}